\documentclass[10pt]{amsart}

\usepackage[english]{babel}

\usepackage{amsmath,amsfonts,amsthm}

\def\ds{\displaystyle}
\def\R{{\mathbb R}}
\def\Z{{\mathbb Z}}

\def\N{{\mathbb N}}

\newtheorem{theo}{Theorem}
\newtheorem{prop}{Proposition}

\newtheorem{lem}{Lemma}

\author[N. Enriquez]{Nathana\"el ENRIQUEZ}
\address{Laboratoire de Probabilit\'es de Paris 6, 4
place Jussieu, 75252 Paris Cedex 05}
\email{enriquez@ccr.jussieu.fr}

\author[C. Sabot]{Christophe Sabot}
\address{
Unit\'e de Math\'ematiques Pures et Appliqu\'ees,
46, All\'ee d'Italie
F-69364 Lyon Cedex 07}
\email{christophe.sabot@umpa.ens-lyon.fr}

\title[Dirichlet environment]{Random walks in a Dirichlet environment}

\begin{document}
\maketitle
{\bf Abstract}:
{\small This paper states a law of large numbers for a random walk in
a random iid environment
on $\Z^d$, where the environment follows some Dirichlet distribution.
Moreover, we give explicit
bounds for the asymptotic velocity of the process and also an
asymptotic expansion of this
velocity at low disorder.}

\section{Introduction}
After a first breakthrough of Kalikow \cite{K},
giving a transience
criterion for  non-reversible multidimensional
Random Walks in Random Environment,
   Sznitman and Zerner proved, several years later,
a law of large numbers in \cite{SZ},
followed by a central limit
theorem proved by Sznitman in \cite{Sznit}. A generalization to the
case of mixing
environments was proved afterwards by Comets and Zeitouni in
\cite{CZ} (we refer to \cite{Zeit} for an overview of the subject).
Despite these progresses,
many important questions, concerning recurrence or explicit criteria
for a ballistic
behavior, remain largely open.

Among random walks in random environment, random walks in an iid
Dirichlet  environment take a
special place, since their annealed law coincides with the law of
some transition reinforced
random walk having an affine reinforcement (see \cite{ES}). These
reinforced walks are defined as follows.
At time 0, we attribute, in a translation invariant way, a weight to
each oriented edge
of $\Z^d$, and  each time the walk crosses an edge, the weight of
this edge is increased by
one. Finally, the walk is a nearest neighbour walk, which chooses, at
each time, an outgoing
edge with a probability which is proportionnal to its weight.

The question of transience and recurrence for such walks, was
answered by Keane and
Rolles, in
\cite{KR}, in the case where the walk evolves on a graph which is a
product of the
integer line with a finite graph. In the context of trees,
a correspondance between reinforced random walks and random
walks in random environment was used before, by Pemantle, in
\cite{Pemantle}. Our purpose is to give some first
results in the
case of $\Z^d$.

In this paper, we state a law of large numbers for such random walks,
under a simple and explicit condition on the weights.
Moreover, we give explicit
bounds for the asymptotic velocity of these walks and also an
asymptotic expansion of this velocity at
low disorder. Low disorder corresponds, in the random environment model,
to the case where the law of the transition probabilities is
concentrated around its mean
value, and, in its reinforcement interpretation, to the case where the
initial weights of the transitions are large, so that these weights are not
significantly affected during the life of the walk (at least, if the
walk is transient).

Let us precise that these walks do not enter the class of walks considered in
\cite{K}, \cite{S} and in several other works, asking the law of the
environment to satisfy a
uniform ellipticity condition.

This ellipticity hypothesis is usually used in two ways :

\_ in the definition of Kalikow's auxiliary Markov chain which involves the
expection of the Green function of the walk killed when exiting a
given set. The
uniform ellipticity is then  a comfortable assumption for checking
the integrability of this
Green function.

\_ in the estimates of the drift of Kalikow's auxiliary Markov chain,
the ellipticity
      condition often plays a key role. We overcome this difficulty by
using an integration by part
formula.

In section 2, we give the definition of random walks in Dirichlet
environment, remind their
connection with transition reinforced random walks and we present our
main results.
In section 3, we present an
integration by part formula that will be the key analytic tool in the
proof of our results.
Indeed, in section 4, it is shown how one can take advantage of the
special form of
the law of the environment, in order to estimate, using the formula
of Section 3, the drift of
the killed Kalikow's auxiliary walk.

In section 5, we study the integrability of the Green
function of the walk which ensures the existence of the original (non
killed) Kalikow's
auxiliary walk and finish the proof of our first result by applying
the law of large
numbers of Sznitman and Zerner \cite{SZ}. In section 6,  we follow
the scheme of
\cite{S} to get bounds for the asymptotic velocity of the walk, and
deduce an expansion of the
asymptotic velocity at low disorder.

\section{Definitions and statement of the results}
We denote by
$T_{2d}:=\{(x_1,...,x_{2d})\in ]0,1]^{2d}, \,
s.t.,\,\sum_{i=1}^{2d}x_i=1\}$, and by  $(e_i)_{1\leq i\leq 2d}$ the
family of unitary vectors of $\Z^d$, defined as follows:
$(e_i)_{1\leq i\leq d}$ is the
canonical basis of $\R^d$, and for all $j\in\{d+1,...,2d\}$ $e_j=-e_{j-d}$.

For all $\vec\alpha:=(\alpha_1,...,\alpha_{2d})\in ]0,+\infty[^{2d}$,
we denote by
$\lambda^{\vec\alpha}$ the Dirichlet probability measure on $T_{2d}$
with parameters
$(\alpha_1,...,\alpha_{2d})$ i.e. the measure on $T_{2d}$:
$${\Gamma(\alpha_1+...+\alpha_{2d})\over
\Gamma(\alpha_1)...\Gamma(\alpha_{2d})}x_1^{\alpha_1-1}...x_{2d}^{\alpha_{2d}-1}dx_1...dx_{2d-1}.$$

For a unit vector $e$ of $\Z^d$, we will sometimes write, for reading
conveniences,
$\alpha_e$ for the weight $\alpha_i$ where $i$ is such that $e_i=e$.

Let us now introduce  random walks in an iid Dirichlet  environment on $\Z^d$.

We define an environment as an element
$\omega=(\omega(x))_{x\in \Z^d}$ where at any
vertex $x$, $\omega(x):=(\omega(x, x+e_1),...,\omega(x, x+e_{2d}))$
belongs to $T_{2d}$.
We set $\mu:=\ds\mathop{\otimes}_{x\in \Z^d}\lambda^{\vec\alpha}$,
so that $\mu$ is a probability measure on the environments such that
$(\omega(x))_{x\in \Z^d}$ are independent
random variables of law $\lambda^{\vec\alpha}$.

We denote by $P^{\omega}$ the law of the Markov chain in the
environment $\omega$ starting at $0$ defined by:
$$\forall x \in \Z^d,\;\; \forall k\in\N,\;\; \forall i=1,...,2d,\;\;
P^{\omega}(X_{k+1}=x+e_i|X_k=x)=\omega(x,x+e_i).$$
The law  of the random walk in random environment (or the so-called
annealed measure) is the
probability measure  $P^\mu=\int P^{\omega}d\mu(\omega)$.

In \cite{ES}, we show that random walks in iid environment have the
law of some reinforced
random walk. The following proposition states that the case of a
Dirichlet environment
corresponds to a quite natural law of reinforcement:

\begin{prop}
The measure  $P^\mu$ satisfies that
$P^\mu$-almost everywhere,
$$P^\mu(X_{n+1}=x+e_i|\sigma
(X_k,k\leq
n))={\alpha_i+N_i(n,X_n)\over\sum_{k=1}^r\alpha_k+N_k(n,X_n)}$$ where
$\vec{N}(n,x)=(N_i(n,x))_{1\leq i\leq 2d}$ and
$N_i(n,x)=\sum_{l=0}^{n-1}1_{\{X_{l+1}-X_l=e_i, X_l=x\}}$.

\end{prop}

We refer the reader to \cite{ES} for the proof.

\subsection{Bounds of the asymptotic velocity}

We can now state our first result:
\begin{theo}
Let $\vec\alpha:=(\alpha_1,...,\alpha_{2d})\in ]0,+\infty[^{2d}$, and
$\mu=\ds\mathop{\otimes}_{x\in
\Z^d}\lambda^{\vec\alpha}$ a probability measure on the environment.
Let us assume that
there exists
$i\in\{1,...,2d\}$ such that $\alpha_{e_i}>1+\alpha_{-e_{i}}$.

The process
$X_n$ is transient under $P^\mu$, and
$$\exists v\in\R^d\setminus\{0\}, \,\mbox{such that}\,P^\mu({X_n\over
n}\to_{n\to\infty} v)=1.
$$
Moreover, for all $i\in\{1,...,d\}$,
$${\alpha_{e_i}-\alpha_{-e_i}-1\over(\sum_{k=1}^{2d} \alpha_k)-1}
\leq
v.e_i\leq {\alpha_{e_i}-\alpha_{-e_i}+1\over(\sum_{k=1}^{2d}
\alpha_k)-1}$$
\end{theo}

\noindent
Remark 1: The assumption on $\alpha_i$ ensures that the set
$$\prod_{i=1}^d[\alpha_{e_i}-\alpha_{-e_i}-1,\alpha_{e_i}-\alpha_{-e_i}+1]$$
does not contain 0. It is a key ingredient in the check of Kalikow's
transience condition.

\noindent Remark 2: When the $\alpha_i$'s are large, $v$ becomes close to the
vector $\ds{1\over\sum_{k=1}^{2d} \alpha_k}
\sum_{i=1}^d(\alpha_{e_i}-\alpha_{-e_i})e_i$. This
is not surprising if one thinks at the corresponding reinforced walk:
the initial weights of the
transitions are large enough so that they are not significantly
affected during the life of the
walk, and the law of the walk becomes close to the law of the Markov chain with
probability transition $\ds{\alpha_i\over\sum_{k=1}^{2d} \alpha_k}$
in the direction
$e_i$.

\noindent
Remark 3: In dimension 1, the condition of theorem 1 is actually
optimal. Indeed from \cite{Sol}, we know that the asymptotic
velocity is not null if and only if either
$E^\mu[{\omega(0,e_1)\over \omega(0, -e_1)}]>1$ or
$E^\mu[{\omega(0,-e_1)\over \omega(0,e_1)}]>1$, which
corresponds exactly to $\alpha_{e_1}>1+\alpha_{-e_1}$ or
$\alpha_{-e_1}>1+\alpha_{e_1}$. Moreover, the asymptotic velocity of
the walk is
equal to ${\alpha_{e_1}-\alpha_{-e_1}-1\over
\alpha_{e_1}+\alpha_{-e_1}-1}$. This
shows the optimality of the lower bound in Theorem 1.

\subsection{Expansion of the velocity in the limit of large parameters}

We turn now to the second result of the paper, which gives the asymptotic
velocity of the walk in
the limit of large parameters $\alpha_k$. Let us remind that, in the
limit of large
parameters $\alpha_k$, the environment is concentrated around its mean value.

Let us fix some notations.
We consider some {\it fixed} transition probabilities
$$ (m_i)\in T_{2d},
$$
and a parameter $\gamma>0$ (aimed to tend to $\infty$).
We consider the weights
$$
\alpha_k=\gamma m_k,
$$
so that the expectation of the transition probability,
$E_\mu(\omega(x,x+e_i)),$
is independent of $\gamma$ and equal to $m_i$.

The mean environment $(m_i)$ defines the transition probabilities of
an homogeneous walk on $\Z^d$, which is ballistic with asymptotic velocity
$$d_m=\sum_{k=1}^{2d} m_{k} e_k,$$
when the mean drift $d_m$ is not null. We denote by $G^m$ its Green function.

The following result gives an estimate in $O({1\over \gamma^2})$
of the asymptotic velocity (in section 6 we give explicit bounds for
this estimate).

\begin{theo}
Assume $d_m\neq 0$.

For $\gamma$ large enough, Theorem 1 applies,
i.e. there exists $v\neq 0$ such that $\lim_{n\to\infty} {X_n\over
n}= v,\, P^\mu\,a.s.$.

Moreover, when $\gamma$ is large,
we have the following  expansion for $v$:
$$v=d_m-{d_m\over \gamma}(G^m(0,0)-1)+O({1\over
\gamma^2}).$$
\end{theo}

\noindent Remark 1: Surprisingly, the second order of the expansion is
colinear to the mean drift $d_m$. We see that $(G^m(0,0)-1)>0$,
which means that there is a slowdown effect, since the second
order term is directed in the opposite direction to the mean drift.

\noindent Remark 2: In \cite{S}, the second author gave an expansion of the
asymptotic velocity in the case of a uniformly elliptic
environment. In this work, several of the estimates relied strongly on the
ellipticity condition, so that the proofs of \cite{S} have here to be
modified. Nevertheless, if we apply the formula of \cite{S} to
this case, we get the same expansion (many simplifications occur
due to the particular expression of the covariance matrix). It is
not surprising that the speed-up effect obtained in some cases of \cite{S}
is not observed in the case of a  Dirichlet environment. The example of
\cite{S}, section 2, was based, indeed, on some correlation between
the transition
probabilities in orthogonal directions. Here, there is a kind of
independence of
the transition probabilities in each direction, in the following
sense: under $\mu=\ds\mathop{\otimes}_{x\in
\Z^d}\lambda^{\vec{\alpha}}$, the law of
$\omega(z,z+e_i)$ is independent of the law of
$({\omega(z,z+e_k)\over 1-\omega(z,z+e_i)})_{k\neq i}$.

\noindent Remark 3: The Green function $G^m(0,0)$ has the following
explicit Fourier
expression
$$G^m(0,0)={1\over(2\pi)^d}\int_{[0,2\pi]^d} {1\over 1-2\sum_{i=1}^d
\sqrt{m_{e_i}m_{-e_i}} \cos(\theta_i) } d\theta_1 \cdots d\theta_d.$$
(we refer to Step 2 of the proof of Proposition 3).
\section{An integration by part formula}

In this section, we present an integration by part formula on
$T_{2d}$ that will appear to be
the key analytic tool in the estimation of the drift of Kalikow's
auxiliary walk.
\begin{lem}
For all $\vec\alpha\in]0,+\infty[^{2d}$, and all differentiable function $f$
on $\R^{2d}$,
$$\int_{T_{2d}}f
d\lambda^{\vec\alpha}={\alpha_1+...+\alpha_d\over\alpha_1}\int_{T_{2d}}x_1.f
d\lambda^{\vec\alpha}+{1\over\alpha_1}\int_{T_{2d}}x_1.((\sum_{k=1}^{2d}
x_k {\partial
f\over\partial x_k}) -{\partial f\over\partial
x_1})d\lambda^{\vec\alpha}.$$
\end{lem}

Proof: We use the well known identity between the Dirichlet law
$\lambda^{\vec\alpha}$
and the
law of the vector $({Z_1\over\sum_{i=1}^{2d} Z_i},
...,{Z_{2d}\over\sum_{i=1}^{2d} Z_i})$  where
the random variables $Z_i$ are independent variables  following the
gamma distribution
$\Gamma(\alpha_i,1)$ of density ${1\over
\Gamma(\alpha_i)}z^{\alpha_i-1}e^{-z}$ on $\R_+$.

This identity implies
$$ \int_{T_{2d}}f d\lambda^{\vec\alpha}= $$
$$
{1\over
\Gamma(\alpha_1)...\Gamma(\alpha_{2d})}
\int_{\R_+^{2d}}f({z_1\over\sum_{i=1}^{2d} z_i},
...,{z_{2d}\over\sum_{i=1}^{2d} z_i})e^{-\sum_{i=1}^{2d} z_i}
z_1^{\alpha_1-1}...z_{2d}^{\alpha_{2d}-1} dz_1...dz_{2d}.$$

Integrating by part with respect to $z_1$, we get
$$ \int_{T_{2d}}f d\lambda^{\vec\alpha}=
{1\over
\Gamma(\alpha_1+1)...\Gamma(\alpha_{2d})}\int_{\R_+^{2d}}(\tilde
f-{\partial \tilde
f\over\partial z_1})e^{-\sum_{i=1}^{2d} z_i}
z_1^{\alpha_1}...z_{2d}^{\alpha_{2d}-1}
dz_1...dz_{2d}$$ where $\tilde
f(z_1,...,z_{2d}):=f({z_1\over\sum_{i=1}^{2d} z_i},
...,{z_{2d}\over\sum_{i=1}^{2d} z_i})$.

Now, we decompose this last integral into  the $\tilde f$-part and
the ${\partial
\tilde f\over\partial z_1}$-part.

Using, in
the reverse sense, the ``Gamma" interpretation of the Dirichlet law
$\lambda^{(\alpha_1+1,\alpha_2,...,\alpha_{2d})}$,
      the $\tilde f$-part becomes
$$ {\Gamma(\alpha_1+...+\alpha_{2d}+1)\over
\Gamma(\alpha_1+1)...\Gamma(\alpha_{2d})}
\int_{T_{2d}}f.x_1^{\alpha_1}...x_{2d}^{\alpha_{2d}-1}dx_1...dx_{2d-1}=
{\alpha_1+...+\alpha_{2d}\over\alpha_1}\int_{T_{2d}}x_1.f
d\lambda^{\vec\alpha}.$$

Now, the  ${\partial
\tilde f\over\partial z_1}$-part writes
$$
-{1\over
\alpha_1\Gamma(\alpha_1)...\Gamma(\alpha_{2d})}\int_{\R_+^d}(z_1.{\partial
\tilde f\over\partial
z_1})e^{-\sum_{i=1}^{2d} z_i}
z_1^{\alpha_1-1}...z_{2d}^{\alpha_{2d}-1} dz_1...dz_{2d}$$
and
$$\begin{array}{rl}\ds z_1.{\partial
\tilde f\over\partial z_1}&\ds =({z_1\over\sum_{i=1}^{2d}
z_i}-{z_1^2\over(\sum_{i=1}^{2d}
z_i)^2})\tilde f_1-{z_1z_2\over(\sum_{i=1}^{2d} z_i)^2}\tilde
f_2-...-{z_1z_{2d}\over(\sum_{i=1}^{2d}
z_i)^2}\tilde f_d\\ &\ds =({z_1\over\sum_{i=1}^{2d} z_i})(\tilde
f_1-\sum_{k=1}^{2d}
({z_k\over\sum_{i=1}^{2d} z_i})\tilde f_k)
\end{array} $$
where $\tilde f_k(z_1,..,z_{2d})=\ds{\partial f\over\partial
x_k}({z_1\over\sum_{i=1}^{2d} z_i},
...,{z_{2d}\over\sum_{i=1}^{2d} z_i})$.

The ``Gamma" interpretation of the Dirichlet law
$\lambda^{(\alpha_1,...,\alpha_{2d})}$ (used for the third time)
allows to conclude.\qed

\section{Kalikow's auxiliary walk}

We remind here the generalization of Kalikow's auxiliary walk (see
\cite{K}) which was already
presented in \cite{S}.

Let $U$ be a connected subset of $\Z^d$, and $\delta\in]0,1]$. We denote by
$\partial U$ the boundary set of $U$, i.e. $\partial U:=
\{z\in\Z^d\setminus U , \exists x\in
U, |z-x|=1\}$.

For all $z\in U$ and $z'\in U\cup\partial U$, and for all environment
$\omega$, we introduce
the Green function of the random walk under the environment $\omega$
killed at rate $\delta$
and at the boundary of $U$:
$$
G_{U,\delta}^\omega(z,z')=E_{z}^{\omega}\left(\sum_{k=0}^{T_U}\delta^k1_{X_k=z'}\right) 
$$
where $T_U=inf\{k, X_k\in\Z^d\setminus U\}$.

{\it In the sequel, we will drop the subscript $\delta$ when
$\delta=1$, and we will write
$G_{U}^\omega(z,z')$ instead of $G_{U,1}^\omega(z,z')$.}

We introduce now the generalized Kalikow's transition probabilities
(originally, Kalikow's transition probabilities were introduced in the
case $\delta=1$):

$$\hat\omega_{U,\delta,z_0}(z,
z+e_i)={E_\mu[G_{U,\delta}^\omega(z_0,z)\omega(z,z+e_i)]\over
E_\mu[G_{U,\delta}^\omega(z_0,z)]}.$$

In order to give bounds for these transition probabilities, we will
be led to apply the
integration by part formula of the previous section to the functions
$G_{U,\delta}^\omega(x,y)$, viewed as functions of the environment $\omega$.

For this purpose, we need the following lemma which gives the expression of the
derivatives of these functions:

\begin{lem}
For all connected subset $U$ of $\Z^d$, for all $x_1,x_2,x_4\in U$,
$x_3\in U\cup \partial U$, $\vert x_3-x_2\vert =1$, and for all
$\delta\in ]0,1[$,
$${\partial G_{U,\delta}^\omega(x_1,x_4)\over\partial (\omega(x_2,x_3))}=
\delta G_{U,\delta}^\omega(x_1,x_2)G_{U,\delta}^\omega(x_3,x_4) $$

\end{lem}
Remark 1: The partial derivative is understood in the following sense:
the function $G_{U, \delta}^\omega$ can be defined by $\sum \delta^n
(\Omega_U)^n$, where $\Omega_U$ is the transition matrix, defined in the proof
below, whose entries are subjected to some stochasticity condition.  But, at
least locally, when $\delta<1$, it can be extended by the same formula to a
function of the variables $(\Omega_U(x,y))$, which are not subjected to
this relation.  In this sense, the partial derivative has a clear
meaning.

Remark 2: When $x_3\in \partial U$, the right-hand term vanishes since
$G^\omega_{U,\delta}(x_3,x_4)=0$.

Proof: Let us define the transition matrix
$\Omega_U(x,y)=\omega(x,y) $ if
$x\in U$, and $\Omega_U(x,y)=0$ if $x\in \partial U$.
We have $$
G_{U,\delta}^\omega(x_1,x_4)=\sum_{n\geq0}\delta^n(\Omega_U)^n_{(x_1,x_4)}$$
and
$${\partial(\Omega_U)^n_{(x_1,x_4)}\over
\partial
(\omega(x_2,x_3))}=\sum_{k_1+k_2=n-1}(\Omega_U)^{k_1}_{(x_1,x_2)}(\Omega_U)^{k_2}_{(x_3,x_4)}$$

so that, taking the derivatives term by term in the sum defining $
G_{U,\delta}^\omega(x_1,x_4)$,
we obtain the result. \qed

We turn now to the estimation of the transition probabilities:
\begin{prop}
For all connected subset $U$ of $\Z^d$, for all $z_0,z\in U$, for all
$\delta\in ]0,1[$ and all $i=1,...,2d$,

$\bullet$ \, if $\ds(\sum_{k=1}^{2d} \alpha_k)>1$, then $\quad
\ds{\alpha_i-1\over(\sum_{k=1}^{2d} \alpha_k)-1}\leq
\hat\omega_{U,\delta,z_0}(z, z+e_i)\leq
{\alpha_i\over(\sum_{k=1}^{2d} \alpha_k)-1}$

$\bullet$ \, if $\ds(\sum_{k=1}^{2d} \alpha_k)<1$, then $\quad
\ds 0\leq
\hat\omega_{U,\delta,z_0}(z, z+e_i)\leq
{\alpha_i-1\over(\sum_{k=1}^{2d}
\alpha_k)-1}$

\end{prop}

Proof: For the clarity of notations we  give the proof for $i=1$.

Lemma 2 yields
$${\partial G_{U,\delta}^\omega(z_0,z)\over\partial (\omega(z, z+e_i))}=
\delta G_{U,\delta}^\omega(z_0,z)G_{U,\delta}^\omega(z+e_i,z) .$$

We now apply Lemma 1 with $f= G_{U,\delta}^\omega(z_0,z)$, viewed as
a function of the
only variables
$x_i:=\omega(z, z+e_i)$ for $i=1,..., 2d$, and we get

$E_\mu[G_{U,\delta}^\omega(z_0,z)]
=
\ds{\alpha_1+...+\alpha_{2d}\over\alpha_1}E_\mu[G_{U,\delta}^\omega(z_0,z)\omega(z,z+e_1)]$

$\hfill{\ds+{1\over\alpha_1}E_\mu\left[\omega(z,
z+e_1).G_{U,\delta}^\omega(z_0,z)\left(\delta\sum_{k=1}^{2d}
\omega(z, z+e_k)G_{U,\delta}^\omega(z+e_k,z)
-\delta G_{U,\delta}^\omega(z+e_1,z)\right)\right]\quad(1) }$

We recall that
$$\delta\sum_{k=1}^{2d}
\omega(z,
z+e_k)G_{U,\delta}^\omega(z+e_k,z)=G_{U,\delta}^\omega(z,z)-1$$
so that the second term in the right side of (1) writes
$${1\over\alpha_1}E_\mu\left[\omega(z,
z+e_1).G_{U,\delta}^\omega(z_0,z)\left(G_{U,\delta}^\omega(z,z)-1
-\delta G_{U,\delta}^\omega(z+e_1,z)\right)\right] $$
so that we get
$$\begin{array}{rl}
E_\mu[G_{U,\delta}^\omega(z_0,z)]=
&\ds{\alpha_1+...+\alpha_{2d}\over\alpha_1}E_\mu[G_{U,\delta}^\omega(z_0,z)\omega(z,z+e_1)]\\
&+\ds{1\over\alpha_1}
E_\mu\left[\omega(z,z+e_1).G_{U,\delta}^\omega(z_0,z)\left(G_{U,\delta}^\omega(z,z)-1
-\delta G_{U,\delta}^\omega(z+e_1,z)\right)\right] \end{array}$$

      and for the ratio $\hat\omega_{U,\delta,z_0}(z,
z+e_1)=\ds{E_\mu[\omega(z, z+e_1)G_{U,\delta}^\omega(z_0,z)]
\over E_\mu[G_{U,\delta}^\omega(z_0,z)]}$,

$\hat\omega_{U,\delta,z_0}(z,
z+e_1)=\ds{\alpha_1\over (\sum_{k=1}^{2d} \alpha_k)-1}$

$\hfill{+\ds{1\over (\sum_{k=1}^{2d} \alpha_k)-1}
{E_\mu\left[\omega(z,
z+e_1).G_{U,\delta}^\omega(z_0,z)\left(G_{U,\delta}^\omega(z,z)
-\delta G_{U,\delta}^\omega(z+e_1,z)\right)\right]\over
E_\mu[G_{U,\delta}^\omega(z_0,z)]}\quad (2)}$

But,
$$\sum_{k=1}^{2d}
\omega(z,
z+e_k)(
G_{U,\delta}^\omega(z,z)-\delta G_{U,\delta}^\omega(z+e_k,z))=1$$
and therefore, for all $k=1,...,2d$,
$$0\leq
G_{U,\delta}^\omega(z,z)-\delta G_{U,\delta}^\omega(z+e_k,z)
\leq{1\over\omega(z,z+e_1)}.$$

These inequalities  allow to bound the ratio in the second term of
the right side of (2),
between  0 and 1, and this finishes the proof.
\qed

\section{Proof of Theorem 1}

We gather now all the ingredients of the proof of Theorem 1. We
want to apply Sznitman
and Zerner's law of large numbers \cite{SZ}. From a careful reading
of the proof of this law of
large numbers, we can see that the only conditions that need to be
fullfilled, are the
integrability of the Green function $G_U^\omega(z_0,z_0)$ for all
bounded $U$, and Kalikow's
condition.

The integrability of the Green function is proved in the following lemma:

\begin{lem}
If there exists $i\in\{1,...,2d\}$, such that $\alpha_i>1$, then
for all connected subset $U$ of $\Z^d$ and all $z_0\in U$,
$E_\mu[G_{U}^\omega(z_0,z_0)]$ is
finite.
\end{lem}

Proof: For the clarity of notations, we suppose that $\alpha_1>1$.

Define now by $N$ the least integer such that $z_0+N e_1$ belongs to
$\partial U$.

We have the following lower bound for the probability $P(\omega, z_0,
U)$ to reach
$\partial U$ from $z_0$ without returning to $z_0$ :
$$P(\omega, z_0, U)\geq \prod_{k=0}^{N-1}\omega(z_0+k e_1,z_0+(k+1)e_1). $$
The number of returns to $z_0$ before hitting $\partial U$, being a
geometric variable whose
parameter is precisely $P(\omega, z_0, U)$, its expectation
$G_{U}^\omega(z_0,z_0)$ is equal to $\ds{1\over P(\omega, z_0, U)}$.

We are now led to examine the integrability of
$E_\mu\left[\left(\prod_{k=0}^{N-1}\omega(z_0+k
e_1,z_0+(k+1)e_1)\right)^{-1}\right]$ which is equal to
$(\int_{T_{2d}}{1\over x_1}
d\lambda^{\vec\alpha})^N$ which is finite since $\alpha_1>1$.\qed

We now have to check Kalikow's  condition.

We notice first that, under the assumption of Theorem 1, Lemma 3
applies and Kalikow's
auxiliary walk is well defined. Then, the monotone convergence
theorem allows to make $\delta$
converge to 1 in the inequalities of Proposition 2.

We then deduce that the drift of Kalikow's
walk belongs to
$$\ds{1\over(\sum_{k=1}^{2d} \alpha_k)-1}
\prod_{i=1}^d[\alpha_{e_i}-\alpha_{-e_i}-1,\alpha_{e_i}-\alpha_{-e_i}+1]$$
which does not contain 0, under the assumption of Theorem 1. This
proves Kalikow's transience
condition.

In order to estimate the asymptotic velocity of the process, we apply
directly Proposition 3.2
of \cite{S} which makes the link between $v$ and the drift of Kalikow's walk.

Remark: in Lemma 3, we only got a sufficient condition for the
integrability of the Green
function to hold. A better result about this question would not have
ameliorated the
statement of Theorem 1 as far as our check of Kalikow's condition
requires a stronger
assumption.

\section{Proof of Theorem 2}

Theorem 2 of section 2 is actually a consequence of a more precise
result, where the
$``O"$ in
the expansion is replaced by an explicit upper bound.

Let us
fix some notations: we set
$$ \gamma=\sum_{i=1}^{2d} \alpha_i,
$$
and
$$
m_i=m_{e_i}={\alpha_i\over \gamma}=E^{\lambda^{(\alpha)}}(\omega (e_i)).
$$
When $\gamma$ is large, the environment $(\omega (x,e_i))$ tends to concentrate
around its mean $(m_i)$,
what can be seen from the expression of the
correlations
$$
\hbox{Cov}_\mu (\omega (x, x+e_i), \omega (x, x+e_j))=\left\{
\begin{array}{l}
- {m_im_j\over \gamma +1} \hbox{, if $i\neq j$}
\\
{m_i (1-\sum_{k\neq i} m_k)\over \gamma +1} \hbox{, if $i= j$},
\end{array}\right.
$$
The mean environment $(m_i)$ defines the transition probabilities of
an homogeneous walk on $\Z^d$, and we define
$$
k_m=
2\sum_{i=1}^d \sqrt{m_{e_i}m_{-e_i}},
$$
so that
$$
1-k_m=\sum_{i=1}^d (\sqrt{m_{e_i}}-\sqrt{m_{-e_i}})^2,
$$
measures the non-symmetry of the walk.
When $k_m<1$, this walk is ballistic with asymptotic velocity
$$
d_m=\sum_{i=1}^{2d} m_{i} e_i,
$$
and we denote by
$G^m(\cdot,\cdot)$ its Green function.
Let us define
$$
\eta_m={{\max_i\sqrt{{m_{e_i}\over m_{-e_i}}}}\over 1-k_m}.
$$

\begin{prop}
Assume we are in the condition of application of Theorem 1, and that
$${2d\over \gamma }\eta_m\le 1,$$
then we have the following estimate
$$\left\vert v-d_m(1-{1\over \gamma-1}(G^m(0,0)-1))\right\vert \le
16\left( {d\over \gamma }\right)^2 {\eta^2_m\over 1-{2d\over
\gamma}\eta_m} .$$

\end{prop}

Proof: Considering the domain
$U=\Z^d$, a killing parameter
$\delta<1$ and $z_0=0$, we get from formula (2)
$$\hat\omega_\delta (z,z+e_i)= m_i+{m_i\over \gamma -1}-
{1\over \gamma -1}
{E_\mu [ G_\delta^\omega (0, z) \omega(z,z+e_i)
(G^\omega_\delta (z,z)-\delta G^\omega_\delta (z+e_i,z)) ]\over E_\mu
[G_\delta (0,z)]}.$$

{\it In the sequel, we will sometimes forget the superscript $\omega$
in $G^\omega_\delta$, when there will be no ambiguity.}

Let us introduce a new probability on the environments $\tilde
\mu(d\omega)$ given by
$$\tilde \mu(d\omega) ={G^\omega_\delta (0,z) \over
E_\mu(G^\omega_\delta(0,z))}\mu(d\omega).$$

We see that
$${E_\mu [G_\delta (0,z)\omega(z,z+e_i)(G_\delta (z,z)-\delta G_\delta
(z+e_i,z))]\over E_\mu [G_\delta (0,z)]}= E_{\tilde\mu}[(G_\delta
(z,z)-\delta G_\delta
(z+e_i,z))\omega(z,z+e_i) ].$$

We proceed as in \cite{S}, and apply Kalikow's formula (cf. the
generalized version in
\cite{S}, Proposition 3.1)  to the measure $\tilde \mu$.

It means that we have
$$E_{\tilde \mu}[ G^\omega_\delta (z,z) \omega(z,z+e_i) ]=
G^{\tilde \omega^z}_\delta (z,z) \tilde \omega^z(z,z+e_i),$$
where $\tilde \omega^z$ is the auxiliary transition probability given by
$$\tilde \omega^z (y,y+e_j)={E_{\tilde\mu} [ G_\delta^\omega
(z,y)\omega(y,y+e_j)]\over E_{\tilde\mu}[
G_\delta^\omega (z,y)]}.$$

Similarly,
$$E_{\tilde \mu}[ G_\delta^\omega (z+e_i,z) \omega(z,z+e_i)]=
G^{\tilde \omega^{z+e_i}}_\delta (z+e_i,z) \tilde \omega^{z+e_i}(z,z+e_i),$$
where $\tilde \omega^{z+e_i}$ is the auxiliary transition probability given by
$$\tilde \omega^{z+e_i} (y,y+e_j)={E_{\tilde\mu}[ G_\delta^\omega
(z+e_i,y)\omega(y,y+e_j)]
\over E_{\tilde \mu}[G_\delta^\omega (z+e_i,y)]}.$$

{\bf Step 1:} We want to estimate the transition probabilities
$\tilde\omega^z$ and
$\tilde \omega^{z+e_i}$.

Lemma 2 yields
$$({\partial\over \partial \omega(y,y+e_k)}- {\partial\over \partial
\omega(y,y+e_j)}) G^\omega_\delta (\cdot ,z)=\delta G^\omega_\delta
(\cdot,y)(G^\omega_\delta
(y+e_k,z)-G^\omega_\delta (y+e_j,z)),$$
moreover
$$\sum_{k=1}^{2d} \omega(y,
y+e_k)(G_\delta^\omega(y,z)-\delta G_\delta^\omega(y+e_k,z))= 1_{y=z}.$$
Using the integration by part formula given
in Lemma 1, we get

$$\begin{array}{rl} m_{e_j} E_{\mu} [G_\delta (0,z) G_\delta (z,y)] =
&E_\mu[ G_\delta
(0,z)G_\delta (z,y)\omega(y,y+e_j)]\\&\\
&\ds+{1 \over \gamma} E_\mu\left[ G_\delta (0,y) \left(G_\delta
(y,z)-\delta G_\delta (y+e_j,z)-1_{y=z}\right)G_\delta
(z,y)\omega(y,y+e_j)\right]\\&\\
& \ds+ {1 \over \gamma} E_\mu\left[ G_\delta (0,z)G_\delta (z,y)
\left(G_\delta (y,y)-\delta G_\delta
(y+e_j,y)-1\right)\omega(y,y+e_j)\right]\end{array}$$

But we have
$$0\le \omega(y,y+e_j)(G_\delta (y,y)-\delta G_\delta (y+e_j,y)) \le
1,\quad(3)$$
and if $y\neq z$
$$\left\vert G_\delta (0,y)\omega(y,y+e_j) (G_\delta (y,z)-\delta
G_\delta (y+e_j,z))\right\vert \le (2d-1) G_\delta (0,z). \quad (4)$$

Indeed, for all $k=1, \cdots, 2d$, we have
$$G_\delta^\omega(y+e_k,z)\ge E^\omega_{y+e_k}[\delta^{T_y}]
G_\delta^\omega(y,z),$$
where $T_y$ is the hitting time of $y$ (equal to infinity if the
random walk never hits $y$).

Since
$${1\over 1-\delta \sum_k \omega(y, y+e_k) E_{y+e_k}^\omega[\delta^{T_y}]}
=G_\delta^\omega(y,y),$$
we get
$$\omega(y,y+e_k)(G_\delta^\omega(y,z)-\delta G_\delta^\omega(y+e_k,z)) \le
{G_\delta^\omega(y,z)\over G_\delta^\omega(y,y)} .$$

But, we also have
$$\sum_{k=1}^{2d} \omega(y,y+e_k)(G_\delta^\omega(y,z)-\delta
G_\delta^\omega(y+e_k,z))=0,
\;\;\; \hbox{if $y\neq z$,}$$
so that we have
$$\left \vert \omega(y,y+e_j)(G_\delta^\omega(y,z)-\delta
G_\delta^\omega(y+e_j,z))
\right\vert \le (2d-1){G_\delta^\omega(y,z)\over G_\delta^\omega(y,y)},$$
which immediately implies the estimate (4).

The inequalities (3) and (4) imply that
$$\left\vert m_{e_j} E_{\mu} [G_\delta (0,z) G_\delta (z,y)] - E_\mu[
G_\delta (0,z)G_\delta (z,y)\omega(y,y+e_j)]\right\vert\\
\leq {2d\over \gamma}E_{\mu} [G_\delta (0,z) G_\delta (z,y)].$$

This gives the following estimate for $\tilde \omega^z$
$$\vert m_{e_j}- \tilde \omega^z(y,y+e_j)\vert \le {2d\over \gamma}.$$

The same procedure gives the same estimate for $\tilde\omega^{z+e_i}$.

Hence, we see that
$$E_{\tilde \mu}[ G_\delta^\omega (z,z) \omega(z,z+e_i) ]=
G^{m+\Delta m}_\delta (z,z) (m_{i}+\Delta m (z,z+e_i)),$$
where $\Delta m(z,z+e_i)$ is a correction to the homogeneous
transition probability $(m_i)$ uniformly bounded by
$$\vert \Delta m\vert  \le {2d\over \gamma}.$$

The same reasoning holds for
$$E_{\tilde \mu}[ G^\omega_\delta (z+e_i,z) \omega(z,z+e_i) ]=
G_\delta^{m+\Delta m}(z+e_i,z) (m_{i}+\Delta m (z,z+e_i)),$$
even if the correction  term $\Delta m$ is not the same.

{\bf Step 2:}
We  compare now the Green function $G_\delta^{m+\Delta m}$ with
$G_\delta^m$.

This is done in \cite{S}, but we reproduce the main lines of the
proof, since we
want to obtain explicit bounds. We first introduce the symmetrizing function
$$\phi^m(z)= \prod_{i=1}^d \sqrt{{m_{e_i}\over m_{-e_i}}}^{z_i}.$$

The Green function $G^m_\delta$ is transformed into
$$G^m_\delta =M_{\phi}^{-1} G_{\delta k_m}^s M_{\phi},\ {(5)}$$
where $M_\phi$ is the operator of multiplication by $\phi$, and
$G^s_{\delta k_m}$ is the Green function of the symmetric random
walk with transition probability
$$s_{e_i}=s_{-e_i}={\sqrt{m_{e_i}m_{-e_i}}\over 2 \sum_{k=1}^{2d}
\sqrt{m_{e_k}m_{-e_k}}}, \;\;\; i=1, \ldots ,d,$$
with killing rates $\delta k_m$ where
$$k_m=2 \sum_{k=1}^{2d}\sqrt{m_{e_k}m_{-e_k}}.$$
We refer to Step 2 of the proof of Lemma 4.3 in \cite{S} for
precisions about this fact.

Hence, we see that
$$\sum_y G_{\delta k}^s(0,y) \le {1\over 1-\delta k_m},$$
which means that
$$\| G_{\delta k}^s \|_{\infty} \le {1\over 1-\delta k_m}.$$

We also have
$$G^{m+\Delta m}_\delta-G^m_\delta = -G^m_\delta
\left(I-(I-\delta \Delta P_m G_\delta^m
)^{-1}\right),$$

where  $\Delta P_m$ is the matrix $(\Delta P_m)_{x,x+e_i}=\Delta m(x,x+e_i)$
(and null anywhere else).

Thus, we get
$$
G^{m+\Delta m}_\delta -G^m_\delta = \delta M_\phi^{-1}
G^s_{\delta k_m}M_\phi \Delta P_m M_\phi^{-1} G^s_{\delta k_m}\left (
I-\delta M_\phi \Delta P_m M_\phi^{-1} G^s_{\delta k_m}\right)^{-1}
M_{\phi},$$

but
$$\| M_\phi \Delta P_m M_\phi^{-1}\|_\infty \le (\max_i \phi(e_i))
{2d\over \gamma}$$
so that we get
$$
(G^{m+\Delta m}_\delta -G^m_\delta)(x,y) \leq \phi(y-x) {2d \over\gamma}
      {1\over 1-k_m} \eta_m
{1\over 1-{2d\over \gamma} \eta_m},$$
and
$$G^{m+\Delta m}_\delta (x,y) \le \phi(y-x) {1\over 1-k_m} {1\over
1-{2d\over \gamma} \eta_m}.$$

This implies that, for all $i\in\{1,...,2d\}$,

$\ds\left\vert \hat\omega_\delta(z,z+e_i) -  m_i + {m_i \over \gamma-1}
(G_\delta^m(0,0)-G_\delta^m(e_i,0)-1)\right\vert$

$\hfill{\begin{array}{l}\leq
\ds{1\over \gamma-1}
\left( 2 {2d\over \gamma} \eta_m^2{1\over 1-{2d\over \gamma}
\eta_m}+ 2 {2d\over \gamma } \eta_m {1\over 1-{2d\over \gamma}\eta_m}
\right)\\
\ds\leq 8 {d\over \gamma^2}{\eta_m^2\over 1-{2d\over
\gamma}\eta_m}\hfill{(6)}\end{array}}$

(we used here  $\eta_m\ge 1$).

The sum
$$
\sum_{i=1}^{2d}
\left(m_i- {m_i \over
\gamma-1}(G_\delta^m(0,0)-G_\delta^m(e_i,0)-1)\right).e_i$$
tends to
$\ds
d_m(1-{1\over \gamma-1}(G^m(0,0)-1))$
when $\delta$ tends to 1. Indeed,
the sum  $\ds\sum_{i=1}^{2d}  m_iG^m(e_i,0).e_i$ cancels,
due to the fact that for each $i\in\{1,...,d\}$
$m_{e_i}G^m(e_i,0)$ and $m_{-e_i}G^m(-e_i,0)$ are both equal to the
common value
$\sqrt{m_{e_i}m_{-e_i}} G^s_k(e_i,0)$ (cf formula (5)).

The triangular inequality combined with the $2d$ inequalities (6)
gives that, for all $z$,
$$
\limsup_{\delta\to 1} \| d_{\hat\omega_\delta}(z)-
d_m(1-{1\over \gamma-1}(G^m(0,0)-1))\|\le
     16 \left( {d\over \gamma}\right)^2 {\eta_m^2\over 1-{2d\over
\gamma}\eta_m},
$$
where $d_{\hat \omega_\delta}(z)= \sum_{k=1}^{2d} \hat\omega(z,z+e_k) e_k$
is the local drift of the transition probability $\hat\omega_\delta$.
     Proposition 3.2 of
\cite{S} allows to conclude.
\qed


\begin{thebibliography}{10}
\bibitem{CZ}
Comets, F., Zeitouni, O.,
A law of large numbers for random walks in random mixing environments.
Ann. Probab.  32  (2004),  no. 1B, 880--914.
\bibitem{ES}
Enriquez, N., Sabot, C.,  Edge oriented reinforced random walks and
RWRE.  C. R.
Math. Acad. Sci. Paris  335  (2002),  no. 11, 941--946.
\bibitem{K}
Kalikow, S., Generalized random walk in a random environment.  Ann.
Probab.  9  (1981),
no. 5, 753--768.
\bibitem{KR}
Keane, M., Rolles, S., Tubular recurrence.  Acta Math. Hungar.  97  (2002),
no. 3, 207--221.
\bibitem{Pemantle}
Pemantle, R., Phase transition in reinforced random walk
and RWRE on trees.  Ann. Probab.  16  (1988),  no. 3,
1229--1241.
\bibitem{S}
Sabot, C., Ballistic random walks in random environment at low
disorder.  Ann. Probab.
32  (2004),  no. 4, 2996--3023.
\bibitem{Sol}
Solomon, F.,
Random walks in a random environment.
Ann. Probability  3  (1975), 1--31.
\bibitem{Sznit}
Sznitman, A-S.,
Slowdown estimates and central limit theorem for random walks in
random environment.
J. Eur. Math. Soc. (JEMS) 2 (2000), no. 2, 93--143.
\bibitem{SZ}
Sznitman, A-S., Zerner, M.,  A law of large numbers for random walks in random
environment.  Ann. Probab.  27  (1999),  no. 4, 1851--1869.
\bibitem{Zeit}
Zeitouni, O.,
Random walks in random environment.  Lectures on probability theory
and statistics,  189--312,
Lecture Notes in Math., 1837, Springer, Berlin, 2004.
\end{thebibliography}
\end{document}